\documentclass[12pt]{amsart}
\usepackage{amssymb,amsmath,amsthm}
\usepackage{comment}
\usepackage[english]{babel}
\usepackage{epsfig}
\usepackage{pstricks}
\usepackage{comment}

\numberwithin{equation}{section}

\newcommand{\Z}{ \mathbb Z}

\newcommand{\C}{ \mathbb C}
\newcommand{\p}[1]{\mathcal{#1}}
\newcommand{\s}[1]{\mathbf{#1}}

\newcommand{\set}[1]{\{#1\}}

\newcommand{\longto}{\longrightarrow}
\DeclareMathOperator{\diag}{diag}
\DeclareMathOperator{\id}{id}
\DeclareMathOperator{\inv}{inv}
\DeclareMathOperator{\cyc}{{cyc}}

\newtheoremstyle{thm}
  {9pt}{9pt}{\itshape}{}{\bfseries}{}{.5em}{}
\theoremstyle{thm}
\newtheorem{thm}{Theorem}[section]
\newtheorem{cor}[thm]{Corollary}
\newtheorem{lemma}[thm]{Lemma}
\newtheorem{prop}[thm]{Proposition}

\newtheoremstyle{defin}
  {9pt}{9pt}{}{}{\bfseries}{}{.5em}{}
\theoremstyle{defin}

\newtheoremstyle{exm}
  {9pt}{9pt}{}{}{\scshape}{}{.5em}{}
\theoremstyle{exm}
\newtheorem{exm}[thm]{Example}

\newtheoremstyle{proof}
  {}{}{}{}{\itshape}{:}{.5em}{}

\title
[A generalization of Foata's fundamental transformation]{A generalization of Foata's fundamental transformation and its applications to the right-quantum algebra}
\author{Matja\v z Konvalinka}
\date{\today}

\begin{document}

\begin{abstract}
 The right-quantum algebra was introduced recently by Garoufalidis, L\^e and Zeilberger in their quantum generalization of the MacMahon master theorem. A combinatorial proof of this identity due to Konvalinka and Pak, and also the recent proof of the right-quantum Sylvester's determinant identity, make heavy use of a bijection related to the first fundamental transformation on words introduced by Foata. This paper makes explicit the connection between this transformation and right-quantum linear algebra identities; applications include a new combinatorial proof of the right-quantum matrix inverse theorem, and two new results, the right-quantum Jacobi ratio theorem and a generalization of the right-quantum MacMahon master thorem.
\end{abstract}

\maketitle

\section{Introduction}

Combinatorial linear algebra is a beautiful and underdeveloped part of enumerative combinatorics.  The underlying idea is very simple: one takes a matrix identity and views it as an algebraic result over a (possibly non-commutative) ring. Once the identity is translated into the language of words, an explicit bijection or an involution is employed to prove the result. The resulting combinatorial proofs are often insightful and lead to extensions and generalizations of the original identities, often in unexpected directions.

\medskip

A tremendous body of literature exists on quantum linear algebra, i.e.\hspace{-0.07cm} on quantum matrices. Without going into definitions, history and technical details let us mention Manin's works \cite{manin-book,manin2}. Recently, the work of Garoufalidis, L\^e and Zeilberger \cite{glz} suggested that certain linear algebra identities (such as the celebrated MacMahon master theorem) are valid in the more general setting of $q$-right-quantum matrices (right-quantum matrices in their terminology). In a series or papers \cite{fh1,fh2,fh3}, Foata and Han reproved the theorem, found interesting further extensions
and an important `$1=q$' principle which allows easy algebraic proofs of certain $q$-equations (implicitly based on the Gr\"obner bases of the underlying quadratic algebras). In a different direction, Hai and Lorenz established the quantum master theorem by using the Koszul duality \cite{hai}, thus suggesting that MacMahon master theorem can be further extended to Koszul quadratic algebras with a large group of (quantum) symmetries.  We refer to \cite{kp} for further references, details and the first combinatorial proof of the right-quantum MacMahon theorem, and some further generalizations. The approach there serves as a basis for \cite{kon}, in which the right-quantum Sylvester's determinant identity is proved by similar means.

\medskip

The main result of this paper (Theorem \ref{main1}) is a non-commutative algebraic identity, whose proof, presented in Section \ref{proof}, is a generalization of combinatorial proofs of crucial arguments in \cite{kp,kon}, and which has numerous applications to right-quantum linear algebra identities. The applications presented are:
\begin{itemize}
 \item the right-quantum matrix inverse formula (Theorem \ref{matinv1}) in Section \ref{matinv},
 \item the right-quantum Jacobi ratio theorem (Theorem \ref{jacobi6}) in Section \ref{jacobi},
 \item a generalization of the right-quantum MacMahon master thorem (Theorem \ref{genmm1}) in Section \ref{genmm}.
\end{itemize}
The method gives new combinatorial proofs of the matrix inverse formula and Jacobi ratio theorem in the commutative case, and Theorem \ref{genmm1} appears to be new even for commutative matrices. We will see that it implies the following Dixon-style identity:
\begin{equation} \label{intro6}
 \sum_{i=1}^{n-1} (-1)^i \binom n{i-1}\binom n{i}\binom n{i+1} = \left\{ \begin{array}{ccl} 2(-1)^m \binom{2m}{m-1}\binom{3m}{m-1} & \colon & n = 2m \\ 0  & \colon & n = 2m-1 \end{array} \right..
\end{equation}

\medskip

The method of proof of Theorem \ref{main1} is related to the \emph{first fundamental transformation} described by Foata in \cite{foata65}.

\section{Notation and the main theorem} \label{main}

Denote by $\p A$ the $\C$-algebra of formal power series with non-commuting variables $a_{ij}$, $1 \leq i,j \leq m$. Elements of $\p A$ are infinite linear combinations of words in variables $a_{ij}$ (with coefficients in $\C$).

\medskip

Words in these variables are often written as \emph{biwords}, i.e.\hspace{-0.07cm} as words in the alphabet $\binom i j$, $1 \leq i,j \leq m$, see for example \cite{fh1}; with this notation, the expression $a_{23}a_{14}a_{22}a_{41}a_{13}$ is written as $\binom{21241}{34213}$. In this paper, however, as in \cite{kp,kon}, we represent such expressions graphically as follows. 

\begin{figure}[ht!]
 \begin{center}
  \includegraphics{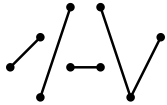}
  \caption{A graphical representation of $a_{23}a_{14}a_{22}a_{41}a_{13}$}
  \label{fig8}
 \end{center}
\end{figure}
We will consider \emph{lattice steps} of the form $(x,i) \to (x+1,j)$ for some $x,i,j \in \Z$, $1 \leq i,j \leq m$.  We think of $x$ being drawn along the $x$-axis, increasing from left to right, and refer to $i$ and $j$ as the \emph{starting height} and \emph{ending height}, respectively. We identify the step $(x,i) \to (x+1,j)$ with the variable $a_{ij}$. Similarly, we identify a finite sequence of steps with a word in the alphabet $\set{a_{ij}}$, $1 \leq i,j \leq m$, i.e.\hspace{-0.07cm} with an element of the algebra $\p A$. Figure \ref{fig8} represents $a_{23}a_{14}a_{22}a_{41}a_{13}$.

\medskip

The \emph{type} of $a_{i_1j_1} a_{i_2j_2} \cdots a_{i_nj_n}$ is defined to be $(\s p;\s r)$ for $\s p = (p_1,\ldots,p_m)$ and $\s r = (r_1,\ldots,r_m)$, where $p_k$ (respectively $r_k$) is the number of $k$'s among $i_1,\ldots,i_n$ (respectively $j_1,\ldots,j_n$). If $\s p  = \s r$, we call the sequence \emph{balanced}.

\medskip

Take non-negative integer vectors $\s p = (p_1,\ldots,p_m)$ and $\s r = (r_1,\ldots,r_m)$ with $\sum p_i = \sum r_i = n$, and a permutation $\pi \in S_m$. An \emph{ordered sequence of type $(\s p;\s r)$ with respect to $\pi$} is a sequence $a_{i_1j_1} a_{i_2j_2} \cdots a_{i_nj_n}$ of type $(\s p;\s r)$ such that $\pi^{-1}(i_k) \leq \pi^{-1} (i_{k+1})$ for $k=1,\ldots,n-1$. Clearly, there are $\binom{r}{r_1,\ldots,r_m}$ elements in $\s O^\pi(\s p;\s r)$, where $r = \sum r_i$. Denote the set of ordered sequence of type $(\s p;\s r)$ with respect to $\pi$ by $\s O^\pi(\s p;\s r)$.

\medskip

A \emph{back-ordered sequence of type $(\s p;\s r)$ with respect to $\pi$} is a sequence $a_{i_1j_1} a_{i_2j_2} \cdots a_{i_nj_n}$ of type $(\s p;\s r)$ such that $\pi^{-1}(j_k) \geq \pi^{-1} (j_{k+1})$ for $k=1,\ldots,n-1$. Denote the set of back-ordered sequences of type $(\s p;\s r)$ with respect to $\pi$ by $\s {\overline O}^\pi(\s p;\s r)$. There are $\binom{p}{p_1,\ldots,p_m}$ elements in $\s {\overline O}^\pi(\s p;\s r)$, where $p = \sum p_i$. 

\begin{exm}
 For $m=3$, $n=4$, $\s p = (2,1,1)$, $\s r = (0,3,1)$ and $\pi = 231$, $\s O^\pi(\s p;\s r)$ is
 $$\set{a_{22}a_{32}a_{12}a_{13},a_{22}a_{32}a_{13}a_{12},a_{22}a_{33}a_{12}a_{12},a_{23}a_{32}a_{12}a_{12}}.$$
 For $m=3$, $n=4$, $\s p = (2,2,0)$, $\s r = (1,2,1)$ and $\pi = 132$, $\s {\overline O}^\pi(\s p;\s r)$ is
 $$\set{a_{12}a_{12}a_{23}a_{21}, a_{12}a_{22}a_{13}a_{21}, a_{12}a_{22}a_{23}a_{11}, a_{22}a_{12}a_{13}a_{21}, a_{22}a_{12}a_{23}a_{11}, a_{22}a_{22}a_{13}a_{11}}.$$
 Figure \ref{fig5} shows some ordered sequences with respect to $1234$ and $2314$, and back-ordered sequences with respect to $1234$ and $4231$.
 
 \begin{figure}[ht!]
  \begin{center}
   \includegraphics{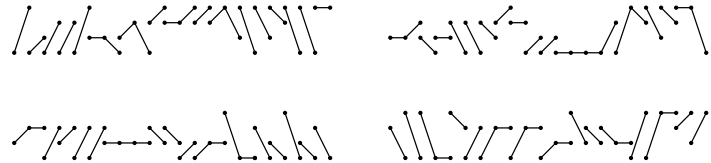}
   \caption{Some ordered and back-ordered sequences.}
   \label{fig5}
  \end{center}
 \end{figure}
\end{exm}

We abbreviate $\s O^\pi(\s p;\s p)$ and $\s {\overline O}^\pi(\s p;\s p)$ to $\s O^\pi(\s p)$ and $\s {\overline O}^\pi(\s p)$, respectively; and if $\pi = \id$, we write simply $\s O(\s p;\s r)$ and $\s {\overline O}(\s p;\s r)$.

\medskip

If each step in a sequence starts at the ending point of the previous step, we call such a sequence a \emph{lattice path}. A lattice path with starting height $i$ and ending height $j$ is called a path from $i$ to $j$.

\medskip

Take non-negative integer vectors $\s p = (p_1,\ldots,p_m)$ and $\s r = (r_1,\ldots,r_m)$ with $\sum p_i = \sum r_i = n$, and a permutation $\pi \in S_m$. Define a \emph{path sequence of type $(\s p;\s r)$ with respect to $\pi$} to be a sequence $a_{i_1j_1} a_{i_2j_2} \cdots a_{i_nj_n}$ of type $(\s p;\s r)$ that is a concatenation of lattice paths with starting heights $i_{k_s}$ and ending heights $j_{l_s}$ so that $\pi^{-1}(i_{k_s}) \leq \pi^{-1} (i_t)$ for all $t \geq k_s$, and $i_t \neq j_{l_s}$ for $t > l_s$. Denote the set of all path sequences of type $(\s p;\s r)$ with respect to $\pi$ by $\s P^\pi(\s p;\s r)$.

\medskip

Similarly, define a \emph{back-path sequence of type $(\s p;\s r)$ with respect to $\pi$} to be a sequence $a_{i_1j_1} a_{i_2j_2} \cdots a_{i_nj_n}$ of type $(\s p;\s r)$ that is a concatenation of lattice paths with starting heights $i_{k_s}$ and ending heights $j_{l_s}$ so that $\pi^{-1}(j_{k_s}) \leq \pi^{-1} (j_t)$ for all $t \leq k_j$, and $j_t \neq i_{k_s}$ for $t < k_s$. Denote the set of all back-path sequences of type $(\s p;\s r)$ by $\s {\overline P}^\pi(\s p;\s r)$. 

\begin{exm}
 Figure \ref{fig4} shows some path sequences with respect to $2341$ and $3421$, and back-path sequences with respect to $1324$ and $4321$. The second path sequence and the second back-path sequence are balanced.
\begin{figure}[ht!]
 \begin{center}
  \includegraphics{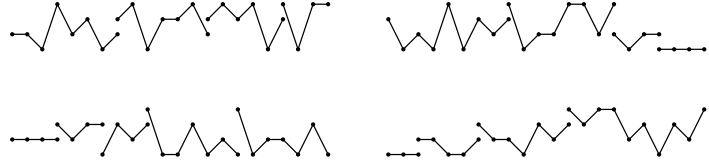}
  \caption{Some path and back-path sequences.}
  \label{fig4}
 \end{center}
 \end{figure}
\end{exm}

We abbreviate $\s P^\pi(\s p;\s p)$ and $\s {\overline P}^\pi(\s p;\s p)$ to $\s P^\pi(\s p)$ and $\s {\overline P}^\pi(\s p)$; and if $\pi = \id$, we write simply $\s P(\s p;\s r)$ and $\s {\overline P}(\s p;\s r)$. Note that a (back-)path sequence of type $(\s p; \s p)$ is a concatenation of lattice paths with the same starting and ending height.

\medskip

For a word $w = i_1i_2\ldots i_n$, say that $(k,l)$ is an inversion of $u$ if $k<l$ and $i_k>i_l$, and write $\inv u$ for the number of inversions of $u$. For $\alpha = a_{i_1j_1} a_{i_2j_2} \cdots a_{i_nj_n}$, write $\inv \alpha = \inv (j_1j_2\ldots j_n) - \inv (i_1i_2\ldots i_n)$. Furthermore, define
$$O^{\pi}(\s p;\s r) = \sum_{\alpha \in \s O^\pi(\s p;\s r)} \alpha, \quad \overline O^{\pi}(\s p;\s r) = \sum_{\alpha \in \s {\overline O}^\pi(\s p;\s r)} (- 1)^{\inv \alpha} \ \alpha,$$
$$P^{\pi}(\s p;\s r) = \sum_{\alpha \in \s P^\pi(\s p;\s r)} \alpha, \quad \overline P^{\pi}(\s p;\s r) = \sum_{\alpha \in \s {\overline P}^\pi(\s p;\s r)} (-1)^{\inv \alpha} \ \alpha,$$

\medskip

Our main theorem seems technical, but it is actually a combinatorial statement with a wide range of applications to right-quantum linear algebra, as we shall see in the following sections. 

\begin{thm} \label{main1}
 Take a matrix $A=(a_{ij})_{m \times m}$, non-negative integer vectors $\s p,\s r$ with $\sum p_i = \sum r_i$, and permutations $\pi,\sigma \in S_m$.
 \begin{enumerate}
  \item Assume that $A$ is \emph{right-quantum}, i.e.\hspace{-0.07cm} that it has the properties
   \begin{eqnarray}
    a_{jk}  a_{ik} & = &  a_{ik}  a_{jk}, \label{main2} \\
    a_{ik}  a_{jl}  -  a_{jk}  a_{il} & = & a_{jl}  a_{ik}  -  a_{il}  a_{jk} \ \ \mbox{for all} \ \ k \neq l. \label{main3}
   \end{eqnarray}
   Then
   \begin{equation} \label{main4}
    O^{\pi}(\s p;\s r) = P^{\sigma}(\s p;\s r).
   \end{equation}
  \item Assume that $A$ satisfies \eqref{main3} above, and that $p_i \leq 1$ for $i=1,\ldots,m$. Then
  \begin{equation} \label{main5}
   \overline O^{\pi}(\s p;\s r) = \overline P^{\sigma}(\s p;\s r).
  \end{equation}
 \end{enumerate}
\end{thm}

\section{Proof of Theorem \ref{main1}} \label{proof}

We can replace $\pi$ by $\id$, since this is just relabeling of the variables $a_{ij}$ according to $\pi$. First we will construct a natural bijection
$$\varphi \colon \s O(\s p;\s r) \longto \s P^\sigma(\s p;\s r).$$

Take an o-sequence $\alpha = a_{i_1j_1} a_{i_2j_2} \cdots a_{i_nj_n}$, and interpret it as a concatenation of steps. Among the steps $i_k \to j_k$ with the lowest $\sigma^{-1}(i_k)$, take the leftmost one. Continue switching this step with the one on the left until it is at the beginning of the sequence. Then take the leftmost step to its right that begins with $j_k$, move it to the left until it is the second step of the sequence, and continue this procedure while possible. Now we have a concatenation of a lattice path and a (shorter) o-sequence. Clearly, continuing this procedure on the remaining o-sequence, we are left with a p-sequence with respect to $\sigma$. 

\begin{exm} The following shows the transformation of
 $$a_{14}a_{12}a_{13}a_{13}a_{14}a_{22}a_{21}a_{23}a_{31}a_{34}a_{33}a_{34}a_{34}a_{34}a_{42}a_{41}a_{42}a_{43}a_{41}a_{41}a_{44}$$
 into 
 $$a_{22}a_{21}a_{14}a_{42}a_{23}a_{31}a_{12}a_{34}a_{41}a_{13}a_{33}a_{34}a_{42}a_{34}a_{43}a_{34}a_{41}a_{13}a_{41}a_{14}a_{44}$$
 with respect to $\sigma = 2341$. In the first five drawings, the step that must be moved to the left is drawn in bold. In the next three drawings, all the steps that will form a path in the p-sequence are drawn in bold.
 \begin{figure}[ht!]
  \begin{center}
   \includegraphics{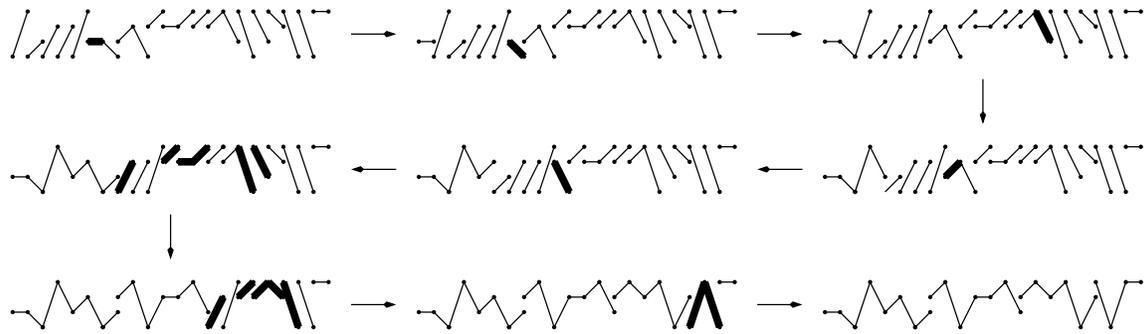}
   \caption{The transformation $\varphi$.}
   \label{fig6}
  \end{center}
 \end{figure}
\end{exm}

\begin{lemma}
 The map $\varphi \colon \s O(\s p;\s r) \to \s P^\sigma(\s p;\s r)$ constructed above is a bijection.
\end{lemma}
\begin{proof}
 Since the above procedure never switches two steps that begin at the same height, there is exactly one o-sequence that maps into a given p-sequence: take all steps starting at height $1$ in the p-sequence
 in the order they appear, then all the steps starting at height $2$ in the p-sequence in the order they appear, etc. Clearly, this map preserves the type of the sequence.
\end{proof}

Define a \emph{q-sequence} to be a sequence we get in the transformation of o-sequences into p-sequences with the above procedure (including the o-sequence and the p-sequence). A sequence $a_{i_1j_1} a_{i_2j_2} \cdots a_{i_nj_n}$ is a q-sequence if it is a concatenation of
\begin{itemize}
 \item some lattice paths with starting heights $i_{k_s}$ and ending heights $j_{l_s}$ so that $\sigma^{-1}(i_{k_s}) \leq \sigma^{-1} (i_t)$ for all $t \geq k_s$, and $i_t \neq j_{l_s}$ for $t > l_s$;
 \item a lattice path with starting height $i_k$ and ending height $j_k$ so that $\sigma^{-1}(i_{k_s}) \leq \sigma^{-1} (i_t)$ for all $t \geq k_s$; and
 \item a sequence that is an o-sequence except that the leftmost step with starting height $j_k$ can be before some of the steps with starting height $i$, $\sigma^{-1}(i) \leq \sigma^{-1}(j_k)$.
\end{itemize} 
For a q-sequence $\alpha$, denote by $\psi(\alpha)$ the q-sequence we get by performing the switch described above; for a p-sequence $\alpha$ (where no more switches are needed), $\psi(\alpha)=\alpha$.  By construction, the map $\psi$ always switches steps that start on different heights.

\medskip

For a sequence $a_{i_1j_1}a_{i_2j_2}\cdots a_{i_nj_n}$, define the \emph{rank} as $\inv (i_1i_2\ldots i_n)$ (more generally, the rank with respect to $\sigma$ is $\inv (\sigma^{-1}(i_1)\sigma^{-1}(i_2)\ldots \sigma^{-1}(i_n))$). Clearly, o-sequences are exactly the sequences of rank $0$. Note also that the map $\psi$ increases by $1$ the rank of sequences that are not p-sequences.

\medskip

Write $\s Q_n^\sigma(\s p;\s r)$ for the union of two sets of sequences of type $(\s p,\s r)$: the set of all q-sequences with rank $n$ and the set of p-sequences (with respect to $\sigma$) with rank $< n$; in particular, $\s O(\s p;\s r)=\s Q_0^\sigma(\s p;\s r)$ and $\s P^\sigma(\s p;\s r)=\s Q_N^\sigma(\s p;\s r)$ for $N$ large enough.

\begin{lemma}
 The map $\psi: \s Q_n^\sigma(\s p;\s r) \to \s Q_{n+1}^\sigma(\s p;\s r)$ is a bijection for all $n$.
\end{lemma}
\begin{proof}
 A q-sequence of rank $n$ which is not a p-sequence is mapped into a q-sequence of rank $n+1$, and $\psi$ is the identity map on p-sequences. This proves that $\psi$ is indeed a map from $\s Q_n^\sigma(\s p;\s r)$ to $\s Q_{n+1}^\sigma(\s p;\s r)$. It is easy to see that $\psi$ is injective and surjective.
\end{proof}

\begin{proof}[Proof of Theorem \ref{main1}]
 Recall that we are assuming that $A$ is right-quantum. Take a q-sequence $\alpha$. If $\alpha$ is a p-sequence, then $\psi(\alpha)=\alpha$. Otherwise, assume that $(x-1,i) \to (x,k)$ and $(x,j) \to (x+1,l)$ are the steps to be switched in order to get $\psi(\alpha)$. If $k = l$, then $\psi(\alpha)=\alpha$ by \eqref{main2}. Otherwise, denote by $\beta$ the sequence we get by replacing these two steps with $(x-1,i) \to (x,l)$ and $(x,j) \to (x+1,k)$. The crucial observation is that $\beta$ is also a q-sequence, and that its rank is equal to the rank of $\alpha$. Furthermore, $\alpha + \beta = \psi(\alpha)+\psi(\beta)$ because of \eqref{main3}. This implies that $\sum \psi(\alpha)=\sum \alpha$ with the sum over all sequences in $\s Q_n^\sigma(\s p;\s r)$. Repeated application of this shows that
 $$\sum \varphi(\alpha) = \sum \alpha$$
 with the sum over all $\alpha \in \s O(\s p;\s r)$. Because $\varphi$ is a bijection, this finishes the proof of \eqref{main4}.\\
 The proof of \eqref{main5} is almost exactly the same. The maps $\psi$ and $\varphi$ must now move steps to the right instead of to the left. Assume that $(x-1,j) \to (x,l)$ and $(x,i) \to (x+1,k)$ are the steps in $\alpha$ we want to switch. The condition $p_i \leq 1$ guarantees that $i \neq j$. Denote by $\beta$ the sequence we get by replacing these two steps with $(x-1,i) \to (x,l)$ and $(x,j) \to (x+1,k)$; $\beta$ is also a q-sequence of the same rank, and because $i \neq j$, its number of inversions differs from $\alpha$ by $\pm 1$. The relation \eqref{main3} implies $\alpha - \beta = \psi(\alpha) - \psi(\beta)$, and this means that $\sum (-1)^{\inv \psi(\alpha)} \psi(\alpha) = \sum (-1)^{\inv \alpha} \alpha$ and hence also
 $$\sum (-1)^{\inv \varphi(\alpha)} \varphi(\alpha) = \sum (-1)^{\inv \alpha} \alpha$$
 with the sum over all $\alpha \in \s {\overline O}(\s p;\s r)$.
\end{proof}

\section{Matrix inverse formula} \label{matinv}

\begin{figure}[ht!]
 \begin{center}
  \includegraphics{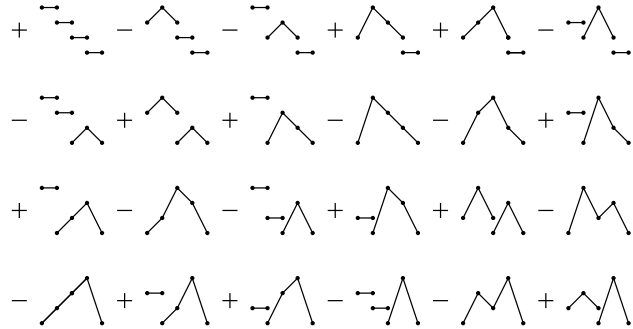}
  \caption{The determinant $\det (a_{ij})_{4 \times 4}$.}
  \label{fig1}
 \end{center}
 \end{figure}

Define the determinant of a matrix $B = (b_{ij})_{m \times m}$ as
$$\det B = \sum_{\pi \in S_m} (-1)^{\inv \pi} b_{\pi(1)1}b_{\pi(2)2}\cdots b_{\pi(m)m}.$$
Note that
$$\det A = \overline O^{w_0}(\s 1),$$
where $A = (a_{ij})_{m \times m}$, $w_0 = m\ldots 21$ and $\s 1 = (1,1,\ldots,1)$. As the first application of Theorem \ref{main1}, we have $\det A = \overline P(\s 1)$ if $A$ is right-quantum; for example, for $m=4$, a graphical representation of $\det A$ is shown in Figure \ref{fig1}.

Note that
\begin{equation} \label{matinv3}
 \det (I-A) = \sum_J (-1)^{|J|} \det A_J,
\end{equation}
where $J$ runs over all subsets of $[m]$ and $A_J$ is the matrix $(a_{ij})_{i,j\in J}$. In other words, $\det (I-A)$ is the weighted sum of $a_{\pi(i_1)i_1}\cdots a_{\pi(i_k)i_k}$ over all permutations $\pi$ of all subsets $\set{i_1,\ldots,i_k}$ of $[m]$, with $a_{\pi(i_1)i_1}\cdots a_{\pi(i_k)i_k}$ weighted by $(-1)^{\cyc \pi}$.

\begin{thm}[right-quantum matrix inverse formula] \label{matinv1}
 If $A=(a_{ij})_{m \times m}$ is a right-quantum matrix, we have
 $$\left(\frac 1 {I-A}\right)_{ij} \, = \, (-1)^{i+j} \cdot \frac 1 {\det(I-A)} \cdot \det\left(I-A\right)^{ji}$$
 for all $i,j$.
\end{thm}

Here $D^{ji}$ means the matrix $D$ without the $j$-th row and $i$-th column.

\medskip

We will prove the equivalent formula
\begin{equation} \label{matinv2}
 \det(I-A) \cdot \left(\frac 1 {I-A}\right)_{ij} \, = \, (-1)^{i+j} \det\left(I-A\right)^{ji}
\end{equation}

If $i=j$, the right-hand side is simply \eqref{matinv3}, with $[m]$ replaced by 
$[m] \setminus \set{i}$, and we can use \eqref{main4} to transform all sequences into p-sequences with respect to $\id$. Figure \ref{fig2} shows the right-hand side of \eqref{matinv2} for $m=4$, $i=j=3$. If $i \neq j$, the right-hand side of \eqref{matinv2} is, again by Theorem \ref{main1}, equal to the sum of all p-sequences with distinct starting and ending heights, with the last lattice path being a path from $i$ to $j$, and with a weight of such a lattice path being $1$ if the number of lattice paths is odd, and $-1$ otherwise. Figure \ref{fig3} shows this for $m=4$, $i=2$, $j=3$. 

\begin{figure}[ht]
 \begin{center}
  \input{rq2.pstex_t}
 \end{center}
 \caption{A representation of $\det \left(I-A\right)^{33}$.}
 \label{fig2}
\end{figure}

\begin{figure}[ht!]
 \begin{center}
  \includegraphics{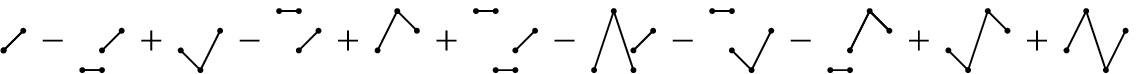}
  \caption{A representation of $- \det \left(I-A\right)^{32}$.}
  \label{fig3}
 \end{center}
\end{figure}

\begin{proof}[Proof of Theorem \ref{matinv1}]
 The left-hand side of \eqref{matinv2} is equal to
 \begin{equation} \label{matinv4}
  \sum (-1)^{\cyc \alpha} \ \alpha \cdot \beta,
 \end{equation}
 where the sum runs over all pairs $(\alpha,\beta)$ with the following properties:
 \begin{itemize}
  \item $\alpha = a_{\pi(i_1)i_1}\cdots a_{\pi(i_k)i_k}$ for some $i_1<\ldots<i_k$, and $\pi$ is a permutation of $\set{i_1,\ldots,i_k}$; $\cyc \alpha$ denotes the number of cycles of $\pi$;
  \item $\beta$ is a lattice path from $i$ to $j$.
 \end{itemize}
 Our goal is to cancel most of the terms and get the right-hand side of \eqref{matinv2}.\\
 Let us divide the pairs $(\alpha,\beta)$ in two groups.
 \begin{itemize}
  \item $(\alpha,\beta) \in \p G_1$ if no starting or ending height is repeated in $\alpha \cdot \beta$, or the first height that is repeated in $\alpha \cdot \beta$ is a starting height;
  \item $(\alpha,\beta) \in \p G_2$ if the first height to be repeated in $\alpha \cdot \beta$ is an ending height.
 \end{itemize}
 The sum \eqref{matinv4} splits into two sums $S_1$ and $S_2$. Let us discuss each of these in turn.
 \begin{enumerate}
  \item Note that if the first height that is repeated in $\alpha \cdot \beta$ is a starting height, this starting height must be $i$, either as the starting height of the first step of $\beta$ if $\alpha$ contains $i$, or the second occurrence of $i$ as a starting height of $\beta$ if $\alpha$ does not contain $i$. For each $\beta$, we can apply \eqref{main5} with respect to $\sigma = (i,1,\ldots,i-1,i+1,\ldots,m)$ to the sum
  $$\sum (-1)^{\cyc \alpha} \alpha$$
  over all $\alpha$ with $(\alpha,\beta) \in \p G_1$. The terms $(-1)^{\cyc \alpha} \ \alpha \cdot \beta$ that do not include $i$ sum up to the right-hand side of \eqref{matinv2}. The terms that do include $i$ either have it in $\alpha$ (and possibly in $\beta$) or they have it only in $\beta$. There is an obvious sign-reversing involution between the former and the latter -- just move the cycle of $\alpha$ containing $i$ over to $\beta$. This means that $S_1$ is equal to the right-hand side of \eqref{matinv2}.
  \item Note that the first height $k$ that is repeated in $\alpha \cdot \beta$ as an ending height cannot be $i$. Fix $k$ and a path $\gamma$ from $k$ to $j$. For each path $\gamma'$ from $i$ to $k$, use \eqref{main5} with respect to $\sigma = (k,1,\ldots,k-1,k+1,\ldots,m)$ on the sum
  $$\sum (-1)^{\cyc \alpha} \ \alpha$$
  over all $\alpha$ such that $(\alpha,\gamma' \gamma) \in \p G_2$ and the only repeated height in $\alpha \cdot \gamma'$ is the ending height $k$. The sum of
  $$\sum (-1)^{\cyc \alpha}\  \alpha \cdot \beta$$
  over $(\alpha,\beta) \in \p G_2$, $\beta=\gamma'\gamma$, and $k$ the only repeated (ending) height in $\alpha \cdot \gamma'$, is therefore equal to
  $$\left(\sum \overline P^\sigma (\s p; \s r)\right) \cdot \gamma$$
  with
  \begin{itemize}
   \item $\s p$ a vector of $1$'s and $0$'s, with $1$ in the $i$-th entry and the $k$-th entry, and
   \item $\s r$ equal to $\s p$ except that the $i$-th entry is $0$ and the $k$-th entry is $2$.
  \end{itemize}
  Equation \eqref{main5} of Theorem \ref{main1} yields
  $$\overline P^\sigma (\s p; \s r) = \overline O(\s p;\s r),$$
  and this is clearly equal to $0$ since $\cdots a_{i'k}a_{i''k} \cdots$ and $\cdots a_{i''k}a_{i'k} \cdots$ have opposite signs in $\overline O(\s p;\s r)$, and since we have \eqref{main2}. This completes the proof.
 \end{enumerate}
\end{proof}

\section{Jacobi ratio theorem} \label{jacobi}

The proof in the previous section is not only the simplest combinatorial proof of the matrix inverse formula (but see \cite{foata79} for an alternative combinatorial proof in the -- less general -- Cartier-Foata case, when $a_{ik}a_{jl} = a_{jl}a_{ik}$ for $i \neq j$), but also generalizes easily to the proof of Jacobi ratio theorem. This result appears to be new (for either Cartier-Foata or right-quantum matrices), although a variant was proved for general non-commutative variables in \cite{gr2} and for quantum matrices in \cite{kl}.

\medskip

We will need the following proposition.

\begin{prop} \label{jacobi1}
 If the matrix $A=(a_{ij})_{m \times m}$ is right-quantum, the matrix $C=(c_{ij})_{m \times m}$ with
 \begin{equation} \label{jacobi7}
  c_{ij} = \left( \frac 1 {I-A}\right)_{ij}
 \end{equation}
 satisfies \eqref{main3}.\\
 Note that $c_{ij}$ is the sum of all paths from $i$ to $j$. 
\end{prop}
\begin{proof}
 We will need some notation:
 \begin{itemize}
  \item let $O$ denote the sum of $O(\s p)$ over all $\s p \geq 0$, and let $P$ denote the sum of $P(\s p)$ over all $\s p \geq 0$;
  \item the superscript $i$ in front of an expression $E$ means that $E$ contains no variable $a_{i*}$; for example, ${}^ic_{ji}$ denotes the sum of all paths from $j$ to $i$ that reach $i$ exactly once, and ${}^{ij}O$ is the sum of $O(\s p)$ over all $\s p \geq 0$ with $p_i = p_j = 0$;
  \item $O_i^j$ for $i \neq j$ means the sum of $O(\s p;\s r)$ with $p_i = r_i + 1$, $p_j = r_j - 1$.
  \item $O_{ij}^k$ for non-equal $i,j,k$ means the sum of $O(\s p;\s r)$ with $p_i = r_i + 1$, $p_j = r_j + 1$, $p_k = r_k - 2$.
  \item $O_{ij}^{kl}$ for non-equal $i,j,k,l$ means the sum of $O(\s p;\s r)$ with $p_i = r_i + 1$, $p_j = r_j + 1$, $p_k = r_k - 1$, $p_l = r_l - 1$.
 \end{itemize} 
 We have to prove $c_{ik}c_{jl} + c_{il}c_{jk} = c_{jk}c_{il} + c_{jl}c_{ik}$ for $k < l$. Let us investigate three possible cases:
 \begin{itemize}
  \item Take $i = k$, $j = l$. First let us prove that $c_{ii} \ {}^ic_{ji} = c_{ji}$. To see this, use \eqref{main4} on $O_i^j$ twice, once with respect to $\pi = (i,j,1,\ldots,i-1,i+1,\ldots,j-1,j+1, \ldots, m)$ and once with respect to $\sigma = (j,i,1,\ldots,i-1,i+1,\ldots,j-1,j+1, \ldots, m)$. We get
  $$O_i^j = c_{ij} \ {}^{j}\! P^\pi \qquad \mbox{and} \qquad O_i^j = c_{jj} \ {}^j c_{ij} \ {}^{j}\! P^\pi$$
  and since ${}^{j}\! P^\pi$ is invertible in $\p A$ (its constant term is $1$), we have
  \begin{equation} \label{jacobi2}
   c_{ij} = c_{jj} \ {}^j c_{ij}.
  \end{equation}
  By using \eqref{main4} on $O$ twice, once with respect to $\pi$ and once with respect to $\sigma$, we get
  \begin{equation} \label{jacobi3}
   O = c_{ii} \ {}^i c_{jj} \ {}^{ij} \! P = c_{jj} \ {}^j c_{ii} \ {}^{ij} \! P.
  \end{equation}
  Furhtermore, \eqref{jacobi2} yields
  $$c_{ii} \ {}^i c_{jj} = c_{ii} \ (c_{jj} - {}^i c_{ji}\ c_{ij}) = c_{ii} \ c_{jj} - (c_{ii} \ {}^i c_{ji}) \ c_{ij} = c_{ii} \ c_{jj} - \ c_{ji}c_{ij}$$
  and similarly
  $$c_{jj} \ {}^j c_{ii} = c_{jj} \ c_{ii} - c_{ij} \ c_{ji},$$
  so \eqref{jacobi3} gives
  $$c_{ii} \ c_{jj} - c_{ji} \ c_{ij} = c_{jj} \ c_{ii} - c_{ij} \ c_{ji}.$$
  \item Take $i \neq k$, $j = l$ (and the case $i = k$, $j \neq l$ is proved analogously). First note that
  $$O_{ij}^k = c_{ik} \ {}^k c_{jk} \ {}^{k}\! P = c_{jk} \ {}^k c_{ik} \ {}^{k} \! P$$
  yields
  \begin{equation} \label{jacobi4}
   c_{ik} \ {}^k c_{jk} = c_{jk} \ {}^k c_{ik}.
  \end{equation}
  Use \eqref{main4} on $O_j^k$ twice, once with respect to $\pi$ and once with respect to $\tau = (j,1,\ldots,j-1,j+1,\ldots,m)$. We get
  $$O_j^k = c_{ii} \ {}^i c_{jk} \ {}^{ik} \! P + c_{ik} \ {}^k c_{ji} \ {}^{ik} \! P \qquad \mbox{and} \qquad O_j^k = c_{jk} \ {}^k \! P = c_{jk} \ {}^k c_{ii} \ {}^{ik} \! P.$$
  But then
  $$c_{ii} \ {}^i c_{jk} + c_{ik} \ {}^k c_{ji} = c_{jk} \ {}^k c_{ii}$$
  implies
  $$c_{ii} \ (c_{jk} - {}^i c_{ji} \ c_{ik}) + c_{ik} \ (c_{ji} - {}^k c_{jk} \ c_{ki}) = c_{jk} \ (c_{ii} - {}^k c_{ik} \ c_{ki})$$
  and
  $$c_{ii} \ c_{jk} + c_{ik} \ c_{ji} = c_{jk} \ c_{ii} + c_{ii} \ {}^i c_{ji} \ c_{ik} + (c_{ik} \ {}^k c_{jk} - c_{jk} \ {}^k c_{ik}) \ c_{ki},$$
  and so \eqref{jacobi2} and \eqref{jacobi4} imply
  $$c_{ii} \ c_{jk} + c_{ik} \ c_{ji} = c_{jk} \ c_{ii} + c_{ji} \ c_{ik}.$$
  \item Assume that $i \neq k$, $j \neq l$. Then
  $$O_{ij}^{kl} = c_{ik} \ {}^k c_{jl} \ {}^{kl} \! P + c_{il} \ {}^l c_{jk} \ {}^{kl} \! P = c_{jk} \ {}^k c_{il} \ {}^{kl} \! P + c_{jl} \ {}^l c_{ik} \ {}^{kl} \!P,$$
  $$c_{ik} \ (c_{jl} - {}^k c_{jk} \ c_{kl}) + c_{il} \ (c_{jk} - {}^l c_{jl} \ c_{lk}) = c_{jk} \ (c_{il} - {}^k c_{ik} \ c_{kl}) + c_{jl} \ (c_{ik} - {}^l c_{il} \ c_{lk}),$$
  and \eqref{jacobi4} imply
  $$c_{ik} \ c_{jl} + c_{il} \ c_{jk} = c_{jk} \ c_{il} + c_{jl} \ c_{ik}.$$
 \end{itemize}
 This completes the proof.
\end{proof}

\begin{thm}[right-quantum Jacobi ratio theorem] \label{jacobi6}
 Take $I,J \subseteq [m]$ with $|I| = |J|$. If $A=(a_{ij})_{m \times m}$ is right-quantum and $C=(c_{ij})_{m \times m}$ is given by
 $$c_{ij} = \left( \frac 1 {I-A}\right)_{ij},$$
 then
 $$\det C_{I,J} = (-1)^{\sum_{i \in I} i + \sum_{j \in J} j} \cdot \frac 1 {\det(I-A)} \cdot  \det(I-A)^{J,I}.$$
 In particular,
 $$\det \left( \frac 1 {I-A}\right) = \frac 1 {\det(I-A)}.$$
\end{thm}

\begin{proof}
 We will only sketch the proof as it is very similar to the proof of Theorem \ref{matinv1} once we have Proposition \ref{jacobi1}, and we will assume that $I=J$ as this makes the reasoning slightly simpler. Use \eqref{main5} (we can do that because of Proposition \eqref{jacobi1}) on $\det C_I$; for a permutation $\pi \in S_m$ with cyclic structure $(i_1^1 i_2^1 \ldots i_{k_1}^1) (i_1^2 i_2^2 \ldots i_{k_2}^2) \cdots (i_1^l i_2^l \ldots i_{k_l}^l)$ (where the first element of each cycle is the smallest, and where starting elements of cycles are increasing), we get the term
 $$(-1)^{\inv \pi} \ (c_{i_1^l,i_2^l}\cdots c_{i_{k_l}^l,i_1^l}) \cdots (c_{i_1^2,i_2^2}\cdots c_{i_{k_2}^2,i_1^2}) \ (c_{i_1^1,i_2^1}\cdots c_{i_{k_1}^1,i_1^1}).$$
 For each selection of paths (in variables $a_{ij}$)
 $$i_1^l \to i_2^l,\ldots,i_{k_l}^l \to i_1^l,\ldots,i_1^2 \to i_2^2,\ldots,i_{k_2}^2 \to i_1^2,\ldots,i_1^1 \to i_2^1,\ldots,i_{k_1}^1 \to i_1^1,$$
 this yields a concatenation of (possibly empty) lattice paths from $i_1^l$ to $i_1^l$, $i_1^{l-1}$ to $i_1^{l-1}$, etc., with exactly one starting height $i_s^t$ marked on each path $i_1^t \to i_1^t$. For example, take $m=5$, $I=\set{1,2,4}$, and $\pi =\binom {124}{421}$. The term of $\det C_I$ corresponding to $\pi$ is $- c_{22} c_{14} c_{41}$, and some of the sequences (without the minus sign) corresponding to this term are depicted in Figure \ref{fig7}. Note the empty path corresponding to $c_{22}$ in the second example.
 \begin{figure}[ht!]
  \begin{center}
   \includegraphics{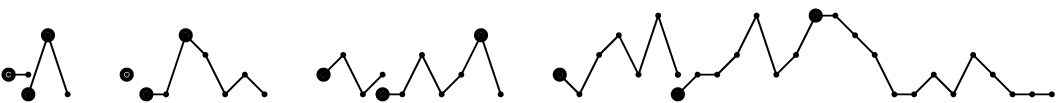}
   \caption{Some sequences in $c_{22} c_{14}c_{41}$.}
   \label{fig7}
  \end{center}
 \end{figure}
 When we multiply $\det C_I$ on the left by $\det(I-A)$, we get a sum
 \begin{equation} \label{jacobi5}
  \sum (-1)^{\cyc \alpha + \inv \beta} \ \alpha \cdot \beta,
 \end{equation}
 where the sum runs over all pairs $(\alpha,\beta)$ with the following properties:
 \begin{itemize}
  \item $\alpha = a_{\pi(i_1)i_1}\cdots a_{\pi(i_k)i_k}$ for some $i_1<\ldots<i_k$, and $\pi$ is a permutation of $\set{i_1,\ldots,i_k}$; $\cyc \alpha$ denotes the number of cycles of $\pi$;
  \item $\beta$ is a concatenation of lattice path from $i_1^l$ to $i_l^1$, $i_1^{l-1}$ to $i_1^{l-1}$ with exactly one starting height $i_s^t$ marked on each path $i_1^t \to i_1^t$, where
  $$\sigma = (i_1^1 i_2^1 \ldots i_{k_1}^1) (i_1^2 i_2^2 \ldots i_{k_2}^2) \cdots (i_1^l i_2^l \ldots i_{k_l}^l)$$
  and $\inv \beta$ denotes the number of inversions of $\sigma$.
 \end{itemize}
 The cancellation process described in the proof of the matrix inverse formula applies here almost verbatim, and this shows that $\det(I-A) \cdot \det C_I$ is equal to $\det (I-A)^I$.
\end{proof}

\section{A generalization of the MacMahon master theorem} \label{genmm}

MacMahon master theorem is a result classically used for proofs of binomial identities. In this section, we will see that the bijection used in the proof of Theorem \ref{main1} gives a far-reaching extension. 

\begin{thm} \label{genmm1}
 Choose a right-quantum matrix $A = (a_{ij})_{m \times m}$, and let $x_1,\ldots,x_m$ be commuting variables that commute with $a_{ij}$. Write $\s x = (x_1,\ldots,x_m)$, $X = \diag \s x$ and $\s x^{\s p}$ for $x_1^{p_1}\cdots x_m^{p_m}$. For $\s p,\s r \geq \s 0$, denote the coefficient
 $$[\s x^{\s r}](A \cdot \s x)^{\s p}$$
 by $G(\s p;\s r)$, and choose an integer vector $\s d$ with $\sum d_i = 0$. Then the generating function
 $$F_A(\s d) = \sum_{\s p = \s r + \s d} G(\s p;\s r) \s x^{\s p}$$
 is a (non-commutative) rational function (with coefficients polynomials in $a_{ij}$) in $x_1,\ldots,x_m$ that can be evaluated as follows. Denote by $\p M$ the multiset of all $i$ with $d_i < 0$, with each $i$ appearing $-d_i$ times, and by $S(\p M)$ the set of all permutations of the multiset $\p M$; denote by $\p N = (N_1,N_2,\ldots,N_d)$ the multiset of all $i$ with $d_i > 0$, with each $i$ appearing $d_i$ times. Note that $\p M$ and $\p N$ have the same cardinality $\delta$; write $M$ (respectively $N$) for the sum of all elements of $\p M$ (respectively $\p N$) counted with their multiplicities. For $\pi=\pi_1\ldots \pi_\delta \in S(\p M)$, let $I_\pi^k$ (for $1 \leq k \leq \delta$) be the set $\set{\pi_1,\ldots,\pi_k}$, write $\varepsilon_\pi^k$ for the size of the intersection of $I_\pi^{k-1}$ and the open interval between $\pi_k$ and $N_k$, and write $J_\pi^k = (I_\pi^k \setminus \set{\pi_k}) \cup \set{N_k}$. Then
 \begin{equation} \label{genmm4}
  F_A(\s d) = \frac{(-1)^{M+N}}{\det(I-XA)} \sum_{\pi \in S(\p M)} \prod_{k=1}^\delta (-1)^{\varepsilon_\pi^k} \cdot  \det(I-XA)^{I_\pi^k,J_\pi^k} \cdot \frac 1 {\det(I-XA)^{I_\pi^k,I_\pi^k}}.
 \end{equation}
\end{thm}

Before proving Theorem \ref{genmm1}, let us show some examples.

\begin{exm}
 If $\s d = \s 0$, $\p M = \p N = \varnothing$, $S(\p M) = \set{\varnothing}$, $M=N=0$, $\delta = 0$ and
 $$F_A(\s 0) = \frac{1}{\det(I-XA)},$$
 which is the right-quantum MacMahon master theorem, see \cite{glz} and \cite{kp}.
\end{exm}

\begin{exm}
 Take
 $$A = \begin{pmatrix} 2 & 1 & 4 & 2 \\ 3 & 2 & 4 & 3 \\ 3 & 4 & 1 & 1 \\ 1 & 3 & 5 & 5 \end{pmatrix},$$
 and $\s d = (1,-2,2,-1)$. i.e.
 $$F(\s d) = 40 xz^2 + 262 x^2z^2 + 128 xyz^2 + 312xz^3 + 251xz^2w + \ldots.$$
 where we write $x,y,z,w$ instead of $x_1,x_2,x_3,x_4$. We have $\p M = (2,2,4)$, $S(\p M)=\set{224,242,422}$, $\p N = (1,3,3)$, $\delta=3$, $M=8$, $N=7$, $I_{224}^1=\set{2}$, $I_{224}^2=\set{2}$, $I_{224}^3=\set{2,4}$, $I_{242}^1=\set{2}$, $I_{242}^2=\set{2,4}$, $I_{242}^3=\set{2,4}$, $I_{422}^1=\set{4}$, $I_{422}^2=\set{2,4}$, $I_{422}^3=\set{2,4}$, $\varepsilon_\pi^i=0$ for all $\pi$ and $i$, $J_{224}^1=\set{1}$, $J_{224}^2=\set{3}$, $J_{224}^3=\set{3,4}$, $J_{242}^1=\set{1}$, $J_{242}^2=\set{2,3}$, $J_{242}^3=\set{3,4}$, $J_{422}^1=\set{1}$, $J_{422}^2=\set{3,4}$, $J_{422}^3=\set{3,4}$.
 Therefore
 \begin{eqnarray*}
  F(\s d) \! = \!  - \frac 1{\det(I-XA)}\biggl( \frac{\det(I-XA)^{2,1}}{\det(I-XA)^{2,2}} \, \frac{\det(I-XA)^{2,3}}{\det(I-XA)^{2,2}} \, \frac{\det(I-XA)^{24,34}}{\det(I-XA)^{24,24}} + \phantom{\biggr) =}\\
   \frac{\det(I-XA)^{2,1}}{\det(I-XA)^{2,2}} \, \frac{\det(I-XA)^{24,23}}{\det(I-XA)^{24,24}} \, \frac{\det(I-XA)^{24,34}}{\det(I-XA)^{24,24}}  + \phantom{xx}\\
   \frac{\det(I-XA)^{4,1}}{\det(I-XA)^{4,4}} \, \frac{\det(I-XA)^{24,34}}{\det(I-XA)^{24,24}} \, \frac{\det(I-XA)^{24,34}}{\det(I-XA)^{24,24}}\phantom{!}\biggr)\phantom{!}=
 \end{eqnarray*}
 $$= - \frac{D_{24,34}\left(D_{2,1}D_{2,3}D_{4,4}D_{24,24} + D_{2,1}D_{2,2}D_{4,4}D_{24,23}D_{24,34} + D_{4,1}D_{2,2}^2D_{24,34}\right)}{D D_{2,2}^2 D_{4,4} D_{24,24}^2},$$
 where
 $$\begin{array}{ccccl}
  D  & = & \det(I-XA) & = & 1 \! - \! 2 x \! - \! 2 y \! - \! z - \! 5 w + \! x y \! - \! 10 x z \! + \! 8 x w- \! 14 y z + \\
  & & & &  y w- \! 5 x y z \! - \! 4 x y w \! + \! 28 x z w \! + \! 17 y z w \! + \! 46 x y z w, \\
  D_{2,1} & = &\det(I-XA)^{2,1} & = & -x - 15 x z - x w + 34 x z w, \\
  D_{2,2} & = &\det(I-XA)^{2,2} & = & 1 - 2 x - z - 5 w - 10 x z + 8 x w + 28 x z w, \\
  D_{2,3} & = & \det(I-XA)^{2,3} & = & -4 z + 5 x z + 17 z w - 30 x z w, \\
  D_{4,1} & = & \det(I-XA)^{4,1} & = & -2 x + x y - 2 x z - 13 x y z, \\
  D_{4,4} & = & \det(I-XA)^{4,4} & = & 1 - 2 x - 2 y - z + x y - 10 x z - 14 y z - 5 x y z, \\
  D_{24,23} & = & \det(I-XA)^{24,23} & = & -z - 4 x z, \\
  D_{24,24} & = & \det(I-XA)^{24,24} & = & 1 - 2 x - z - 10 x z, \\
  D_{24,34} & = & \det(I-XA)^{24,34} & = & -4 z + 5 x z.
 \end{array}$$
 \end{exm}
 
\begin{proof}[Proof of Theorem \ref{genmm1}]
 Fix non-negative integer vectors $\s p,\s r$ with $\s p = \s r + \s d$, and use \eqref{main4} on $\s O(\s p;\s r)$ with respect to the permutation $\sigma=i_1\cdots i_sj_1\cdots,j_t$, where $i_1<i_2<\ldots<i_s$ form the underlying set of $\p N$ (in other words, $\set{i_1,\ldots,i_s} = \set{i \colon d_i>0}$) and $j_1<j_2<\ldots<j_t$ are the remaining elements of $\set{1,\ldots,m}$. A path sequence in $\s P^\sigma(\s p;\s r)$ has the following structure. The first path starts at $N_1=i_1$ and ends at one of the heights in $\p M$; the second path starts at $N_2$ (which is $i_1$ if $d_{i_1} > 1$, and $i_2$ if $d_{i_1} = 1$), and ends at one of the heights in $\p M$, and it does not include the ending height of the previous path except possibly as the ending height. In general, the $k$-th path starts at $N_k$ and ends at one of the heights in $\p M$, and does not contain any of the ending heights of previous paths except possibly as the ending height. All together, the ending heights of these $\delta$ paths form a permutation of $\p M$, which explains why $F_A(\s d)$ is written as a sum over $\pi \in S(\p M)$. After these paths, we have a balanced path sequence that does not include any height in $\p M$.\\
 Now choose $\pi = \pi_1\cdots \pi_\delta \in S(\p M)$, and look at all the p-sequences in $\s P^\sigma(\s p;\s r)$ (for all $\s p,\s r \geq \s 0$ with $\s p = \s r + \s d$) whose first $\delta$ ending heights of paths are $\pi_1,\ldots,\pi_\delta$ (in this order). The $k$-th path is a path from $N_k$ to $\pi_k$, and it does not include $\pi_1,\ldots,\pi_{k-1}$ except possibly as an ending height. By the inverse matrix formula, such paths, weighted by $\s x^{\s {p_k}}$ where $(\s {p_k},\s {r_k})$ is the type of the path, are enumerated by
 $$ \pm \frac{1}{\det(I-XA)^{I_\pi^{k-1},I_\pi^{k-1}}} \cdot \det(I-XA)^{I_\pi^k,J_\pi^k},$$
 and a simple consideration shows that the sign is $(-1)^{N_k+\pi_k+\varepsilon_\pi^k}$.  The balanced path sequences that do not include heights from $\p M$ are enumerated by
 $$\!\left({\textstyle \frac {1} {\det(I-XA)^{I_\pi^\delta}} \! \cdot \! \det(I-XA)^{I_\pi^\delta \cup \set{j_1}}} \right) \cdot \left( {\textstyle \frac 1{\det(I-A)^{I_\pi^\delta \cup \set{j_1}}}} \cdot \det\left(I-A\right)^{I_\pi^\delta \cup \set{j_1,j_2}}\right) \cdots \frac 1{1-a_{j_tj_t}} =$$
 $$= \frac 1 {\det(I-XA)^{I_\pi^\delta}},$$
 where we wrote $\det(I-XA)^I$ instead of $\det(I-XA)^{I,I}$. Formula \eqref{genmm4} follows.
\end{proof}

\begin{proof}[Proof of \eqref{intro6}]
 Let us denote the sum we are trying to calculate by $S(n)$. Clearly,
 $$[x^{n+1} y^{n+1} z^{n-2}] (z-y)^n(x-z)^n(y-x)^n = [xyz^{-2}] ( 1 - {\textstyle \frac y z})^n( 1 - {\textstyle \frac z x})^n( 1 - {\textstyle \frac x y})^n = $$
 $$= [xyz^{-2}] \sum_{i,j,k} (-1)^{i+j+k} \binom n i\binom n j\binom n k \left(\frac x y \right)^i \left(\frac y z \right)^j \left(\frac z x \right)^k =S(n),$$
 and so we have to use Theorem \ref{genmm1} for 
 $$A = \begin{pmatrix} 0 & -1 & 1 \\ 1 & 0 & -1 \\ -1 & 1 & 0 \end{pmatrix}$$
 and $\s d = (-1,-1,2)$. We get
 $$\sum_{\s p = \s r + \s d} G(\s p;\s r) x^{p_1} y^{p_2} z^{p_3} = \frac{-z^2 (1+x)}{(1+x z) (1+xy+xz+yz)}+\frac{-z^2 (1-y)}{(1+y z) (1+ xy+xz+yz)}$$
 and
 $$S(n) = [x^n y^n z^n] \left( \frac{-z^2 (1+x)}{(1+x z) (1+xy+xz+yz)}+\frac{-z^2 (1-y)}{(1+y z) (1+ xy+xz+yz)} \right) =$$
 $$= 2 [x^n y^n z^n] \left( \frac{-z^2}{(1+x z) (1+xy+xz+yz)} \right),$$
 where we used some obvious symmetry. Then
 $$S(n) = 2 \sum_{i,j,k,l} (-1)^{l+1} (xz)^l z^2 \binom{i+j+k}{i,j,k} (-1)^{i+j+k} (xy)^i(xz)^j(yz)^k,$$
 with the sum over all $i,j,k,l \geq 0$ with $l+i+j = n$, $i+j = n$, $l+2+j+k = n$, i.e.\hspace{-0.07cm} $S(n) = 0$ if $n$ is odd and
 $$S(2m) = \frac{2(-1)^m}{(m+1)!(m-1)!}  \sum_{l=0}^{m-1} \frac{(3m-1-l)!}{(m-1-l)!} = 2 (-1)^m \binom{2m}{m-1} \sum_{l=0}^{m-1} \binom{3m-1-l}{m-1-l},$$
 and so $S(2m) = 2(-1)^m \binom{2m}{m-1}\binom{3m}{m-1}$ since every $(m-1)$-subset of $\set{1,\ldots,3m}$ consists of elements $\set{1,\ldots,l}$ and an $(m-1-l)$-subset of $\set{l+2,\ldots,3m}$ for a uniquely determined $l$.
\end{proof}

\section{Final remarks}

Some of the results (Theorems \ref{main1} and \ref{matinv1}) have natural $q$- and $\s q$-analogues; these can either be proved by the $1=q$ and $1 = q_{ij}$ principles (see \cite[\S 3]{fh1} and \cite[Lemma 12.4]{kp}) or by some straightforward bookkeeping, cf.\hspace{-0.07cm} \cite[\S\S 5--8]{kp}. However, Theorems \ref{jacobi6} and \ref{genmm1} do not seem to extend to a formula for general $q$ or $q_{ij}$.

\medskip

As can be seen from the proof of Theorem \ref{matinv1}, the formula
$$\left(\frac 1 {I-A}\right)_{mm} \, = \, \frac 1 {\det(I-A)} \, \cdot \, \det\left(I-A\right)^{mm}$$
holds even if \eqref{main2} and \eqref{main3} hold only for $k,l \neq m$. This fact (in the more general $\s q$-right-quantum case) was used in the proof of $\s q$-Cartier-Foata and $\s q$-right-quantum Sylvester's determinant identity \cite[Theorem 1.3]{kon}. 

\medskip

It is easy to prove that the variables $c_{ij}$ defined by \eqref{jacobi7} satisfy \eqref{main2}, i.e.\hspace{-0.07cm} that the matrix $C=(c_{ij})_{m \times m}$ is right-quantum. We do not need this fact, however.

\medskip

Even though the statement of Theorem \ref{genmm1} appears rather intricate even in the commutative case, a computer algebra program finds the generating function easily. A Mathematica package {\tt genmacmahon.m} that calculates $f_A(\s d)$ for a commutative matrix $A = (a_{ij})_{m \times m}$ and an integer vector $\s d = (d_1,\ldots,d_m)$ with $\sum d_i = 0$ is available at {\tt{http://www-math.mit.edu/\~{}konvalinka/genmacmahon.m}} (read in the package with {\tt << genmacmahon.m} and write {\tt F[A,d,x]} to calculate the generating function $f_A(\s d)$ in variables $x_1,\ldots,x_m$), and it would be easy to adapt this to the non-commutative situation.

\medskip

It would be nice to use Thoerem \ref{genmm1} to prove
\begin{equation} \label{genmm2}
 \sum_{i=k}^{n-k} (-1)^i \binom n{i-k}\binom n{i}\binom n{i+k} = \frac{(-1)^m(2m)!^2(3m)!}{m!(m - k)!(m + k)!(2m - k)!(2m + k)!},
\end{equation}
which can be established by the WZ method \cite{pwz}. The author proved that
$$S_k(n)=\sum_{i=k}^{n-k} (-1)^i \binom n{i-k}\binom n{i}\binom n{i+k}$$
satisfies $S_k(n)=0$ if $n$ is odd and, if $k \geq 1$,
\begin{equation} \label{genmm3}
 S_k(2m) = 2 \sum_{j=1}^k \! \binom{2k-j-1}{k-1} \! \sum_{i=0}^{\lfloor j/2 \rfloor} \! (-1)^{m-i} \binom j{2i} \! \binom{3m-i+j-k}{m-i} \! \binom{2m}{m+k-i},
\end{equation}
and it is easy to see for small $k$ that this is equal to the right-hand side of \eqref{genmm2}.

\bigskip

{\bf Acknowledgment.} The author is grateful to Igor Pak for his help, suggestions and support.

\bibliography{rq}
\bibliographystyle{amsalpha}

\bigskip

\bigskip

\bigskip

{\sc \scriptsize Department of Mathematics, Massachusetts Institute of Technology, Cambridge, MA 02139\\
\tt{konvalinka@math.mit.edu}\\
\tt{http://www-math.mit.edu/\~{}konvalinka/}}

\end{document}